\documentclass[a4paper,12pt]{article}
    \usepackage[top=2.5cm,bottom=2.5cm,left=2.5cm,right=2.5cm]{geometry}
    \usepackage{cite, amsmath, amssymb}
    \usepackage[margin=1cm,%
                font=small,%
                format=hang,%
                labelsep=period,%
                labelfont=bf]{caption}
\usepackage[latin1]{inputenc}
\usepackage{rotating}
\usepackage{color}
\usepackage{graphicx}
\usepackage{epstopdf}
\usepackage{multirow} % para las tablas
\usepackage{subfigure}
\usepackage{amsmath}
\usepackage{amsfonts}
\usepackage{latexsym}
\usepackage{authblk}
\usepackage[none]{hyphenat}
\begin{document}

\sloppy
\newcommand{\hide}[1]{}
\newcommand{\tbox}[1]{\mbox{\tiny #1}}
\newcommand{\half}{\mbox{\small $\frac{1}{2}$}}
\newcommand{\sinc}{\mbox{sinc}}
\newcommand{\const}{\mbox{const}}
\newcommand{\trc}{\mbox{trace}}
\newcommand{\intt}{\int\!\!\!\!\int }
\newcommand{\ointt}{\int\!\!\!\!\int\!\!\!\!\!\circ\ }
\newcommand{\eexp}{\mbox{e}^}
\newcommand{\bra}{\left\langle}
\newcommand{\ket}{\right\rangle}
\newcommand{\EPS} {\mbox{\LARGE $\epsilon$}}
\newcommand{\ar}{\mathsf r}
\newcommand{\im}{\mbox{Im}}
\newcommand{\re}{\mbox{Re}}
\newcommand{\bmsf}[1]{\bm{\mathsf{#1}}}
\newcommand{\mpg}[2][1.0\hsize]{\begin{minipage}[b]{#1}{#2}\end{minipage}}
\newcommand{\CC}{\mathbb{C}}
\newcommand{\NN}{\mathbb{N}}
\newcommand{\PP}{\mathbb{P}}
\newcommand{\RR}{\mathbb{R}}
\newcommand{\QQ}{\mathbb{Q}}
\newcommand{\ZZ}{\mathbb{Z}}
%%%%%%%%%%%%%%%%%%%
% Letras griegas  %
%%%%%%%%%%%%%%%%%%%
\renewcommand{\a}{\alpha}
\renewcommand{\b}{\beta}
\renewcommand{\d}{\delta}
\newcommand{\D}{\Delta}
\newcommand{\g}{\gamma}
\newcommand{\G}{\Gamma}
\renewcommand{\th}{\theta}
\renewcommand{\l}{\lambda}
\renewcommand{\L}{\Lambda}
\renewcommand{\O}{\Omega}
\newcommand{\s}{\sigma}
%%%%%%%%%%%%%

\newtheorem{theorem}{Theorem}
\newtheorem{acknowledgement}[theorem]{Acknowledgement}
\newtheorem{algorithm}[theorem]{Algorithm}
\newtheorem{axiom}[theorem]{Axiom}
\newtheorem{claim}[theorem]{Claim}
\newtheorem{conclusion}[theorem]{Conclusion}
\newtheorem{condition}[theorem]{Condition}
\newtheorem{conjecture}[theorem]{Conjecture}
\newtheorem{corollary}[theorem]{Corollary}
\newtheorem{criterion}[theorem]{Criterion}
\newtheorem{definition}[theorem]{Definition}
\newtheorem{example}[theorem]{Example}
\newtheorem{exercise}[theorem]{Exercise}
\newtheorem{lemma}[theorem]{Lemma}
\newtheorem{notation}[theorem]{Notation}
\newtheorem{problem}[theorem]{Problem}
\newtheorem{proposition}[theorem]{Proposition}
\newtheorem{remark}[theorem]{Remark}
\newtheorem{solution}[theorem]{Solution}
\newtheorem{summary}[theorem]{Summary}
\newenvironment{proof}[1][Proof]{\noindent\textbf{#1.} }{\ \rule{0.5em}{0.5em}}
%%%%%%%%%%%%%%%%%%%%%%%%%%%%%%%%%%%%%%%%%%%5

\title{\vspace*{3.5cm}\textbf{Revan--degree indices on random graphs}}

\author[1]{\bf{R. Aguilar-S\'anchez}} 
\author[2]{\bf{I. F. Herrera-Gonz\'alez}}
\author[3]{\bf{J. A. M\'endez-Berm\'udez\footnote{Corresponding author}}}
\author[4]{\bf{Jos\'e M. Sigarreta}}

\affil[1]{\it Facultad de Ciencias Qu\'imicas, Benem\'erita Universidad Aut\'onoma de Puebla,
Puebla 72570, Mexico}
\affil[2]{\it Departamento de Ingenier\'ia, Universidad Popular Aut\'onoma del 
Estado de Puebla, Puebla, Pue., 72410, Mexico}
\affil[3]{\it Instituto de F\'{\i}sica, Benem\'erita Universidad Aut\'onoma de Puebla,
Apartado Postal J-48, Puebla 72570, Mexico}
\affil[4]{\it Universidad Aut\'onoma de Guerrero, Centro Acapulco CP 39610,
Acapulco de Ju\'arez, Guerrero, Mexico}
\affil[ ]{\ttfamily {\textbf {ras747698@gmail.com, ivanfernando.herrera@upaep.mx, jmendezb@ifuap.buap.mx, josemariasigarretaalmira@hotmail.com}}}
\date{}
\maketitle
\thispagestyle{empty}

\centerline{(Received June 5, 2020)}
\begin{abstract}
Given a simple connected non-directed graph $G=(V(G),E(G))$, we consider two families of graph 
invariants: 
$RX_\Sigma(G) = \sum_{uv \in E(G)} F(r_u,r_v)$ (which has gained interest recently) and 
$RX_\Pi(G) = \prod_{uv \in E(G)} F(r_u,r_v)$ (that we introduce in this work);
where $uv$ denotes the edge of $G$ connecting the vertices $u$ and $v$, $r_u$ is the Revan  
degree of the vertex $u$, and $F$ is a function of the Revan vertex degrees. Here,
$r_u = \Delta + \delta - d_u$ with $\Delta$ and $\delta$ the maximum and minimum degrees 
among the vertices of $G$ and $d_u$ is the degree of the vertex $u$.
Particularly, we apply both $RX_\Sigma(G)$ and R$X_\Pi(G)$ on two models of random graphs: 
Erd\"os-R\'enyi graphs and random geometric graphs.
By a thorough computational study we show that $\left< RX_\Sigma(G) \right>$ and 
$\left< \ln RX_\Pi(G) \right>$, normalized to the order of the graph, scale with the average Revan 
degree $\left< r \right>$; here $\left< \cdot \right>$ denotes the average over an ensemble of random 
graphs.
Moreover, we provide analytical expressions for several graph invariants of both families in the 
dense graph limit.
\end{abstract}

\baselineskip=0.30in
\section{Introduction}

We can identify two families of graph invariants which have been extensively studied in chemical 
graph theory, namely  
\begin{equation}
X_\Sigma(G) = \sum_{uv \in E(G)} F(d_u,d_v)
\label{TI}
\end{equation}
and 
\begin{equation}
X_\Pi(G) = \prod_{uv \in E(G)} F(d_u,d_v) \, .
\label{MTI}
\end{equation}
Here $uv$ denotes the edge of the graph $G=(V(G),E(G))$ connecting the vertices $u$ and $v$,
$d_u$ is the degree of the vertex $u$, and $F(x,y)$ is a given function of the vertex degrees, 
see e.g.~\cite{G13}. 
While both $X_\Sigma(G)$ and $X_\Pi(G)$ are referred as topological indices in the literature, to make
a distinction between them, here we name $X_\Sigma(G)$ and $X_\Pi(G)$ as topological indices (TIs) 
and multiplicative topological indices (MTIs), respectively.

In fact, within a statistical approach on random graphs, it has been recently shown that the average 
values of indices of the type $X_\Sigma(G)$, normalized to the order of the graph $n$, scale with the 
average degree $\left< d \right>$; see e.g. Refs.~\cite{MMRS20,AHMS20,MMRS21,AMRS21,MAAS22}. 
That is, $\left< X_\Sigma(G) \right>/n$ is a function of $\left< d \right>$ only:
\begin{equation}
\frac{\left< X_\Sigma(G) \right>}{n} \equiv f_\Sigma(\left< d \right>) \, .
\label{scaling1}
\end{equation}

More recently, a number of new TIs with the form
\begin{equation}
RX_\Sigma(G) = \sum_{uv \in E(G)} F(r_u,r_v)
\label{RTI}
\end{equation}
have been proposed and studied, see e.g.~Refs.~\cite{K17,K18,KG22,KMRS22}. 
Above, $r_u$ is the {\it Revan vertex degree} of the vertex $u$ which is defined as
\begin{equation}
r_u = \Delta + \delta - d_u ,
\label{ru}
\end{equation}
where $\Delta$ and $\delta$ are the maximum and minimum degrees among the vertices of the graph $G$,
respectively. Note that $RX_\Sigma(G)$ is the Revan version of $X_\Sigma(G)$.

Thus, inspired by the scaling law of Eq.~(\ref{scaling1}), in this paper we explore the statistical properties 
of $\left< RX_\Sigma(G) \right>$ on random graphs and look for the scaling parameter and the 
corresponding scaling law.
Moreover, to complete the panorama of Revan-degree--based indices, we introduce Revan versions of MTIs:
\begin{equation}
RX_\Pi(G) = \prod_{uv \in E(G)} F(r_u,r_v) \, ,
\label{RMTI}
\end{equation}
and also study their statistical and scaling properties on random graphs.

\section{Statistical analysis of Revan-degree--based TIs on random graphs}
\label{Sec:TI}

Among the recently introduced Revan-degree--based indices, $RX_\Sigma(G)$, we can mention~\cite{K17,K18,KG22}
$$
R_1(G) = \sum_{uv\in E(G)} r_u + r_v \, ,
\qquad
R_2(G) = \sum_{uv\in E(G)} r_u r_v \, ,
$$
$$
FR(G) = \sum_{uv\in E(G)} r_u^2 + r_v^2 \, ,
$$
and
$$
RSO(G) = \sum_{uv\in E(G)} \sqrt{r_u^2 + r_v^2} \, .
$$
Evidently, these TIs are the Revan versions of the first and second Zagreb indices~\cite{GT72},
$$
M_1(G) = \sum_{uv\in E(G)} d_u + d_v \, ,
\qquad
M_2(G) = \sum_{uv\in E(G)} d_u d_v \, ,
$$
the forgotten index~\cite{FG15}
$$
F(G) = \sum_{uv\in E(G)} d_u^2 + d_v^2 \, ,
$$
and the Sombor index~\cite{G21}
$$
SO(G) = \sum_{uv\in E(G)} \sqrt{d_u^2 + d_v^2} \, ,
$$
respectively.

In what follows we compute $R_1(G)$, $R_2(G)$, $FR(G)$ and $RSO(G)$ on two models of random graphs: 
Erd\"os-R\'enyi (ER) graphs and random geometric (RG) graphs. 
ER graphs~\cite{SR51,ER59} $G_{\tbox{ER}}(n,p)$ are formed by $n$ vertices connected independently 
with probability $p \in [0,1]$. 
While RG graphs~\cite{DC02,P03} $G_{\tbox{RG}}(n,r)$ consist of $n$ vertices uniformly and independently 
distributed on the unit square, where two vertices are connected by an edge if their Euclidean distance is less 
or equal than the connection radius $\ell \in [0,\sqrt{2}]$.
Moreover, since a given parameter pair [$(n,p)$ or $(n,\ell)$] represents an infinite-size ensemble of random 
[ER or RG] graphs, the computation of a graph invariant on a single graph may be irrelevant. In contrast, 
the computation of the average value of a graph invariant over a large ensemble of random graphs, all 
characterized by the same parameter pair, may provide useful 
{\it average} information about the full ensemble. This {\it statistical} approach, well known in random matrix 
theory studies, has been recently applied to random graphs and networks by means of degree--based TIs, 
see e.g. Refs.~\cite{MMRS20,AHMS20,MMRS21,AMRS21,MAAS22,KMRS22}.

\subsection{Revan-degree--based TIs on Erd\H{o}s-R\'{e}nyi graphs}
\label{Sec:TIonER}

\begin{figure}
\centering
\includegraphics[width=0.8\textwidth]{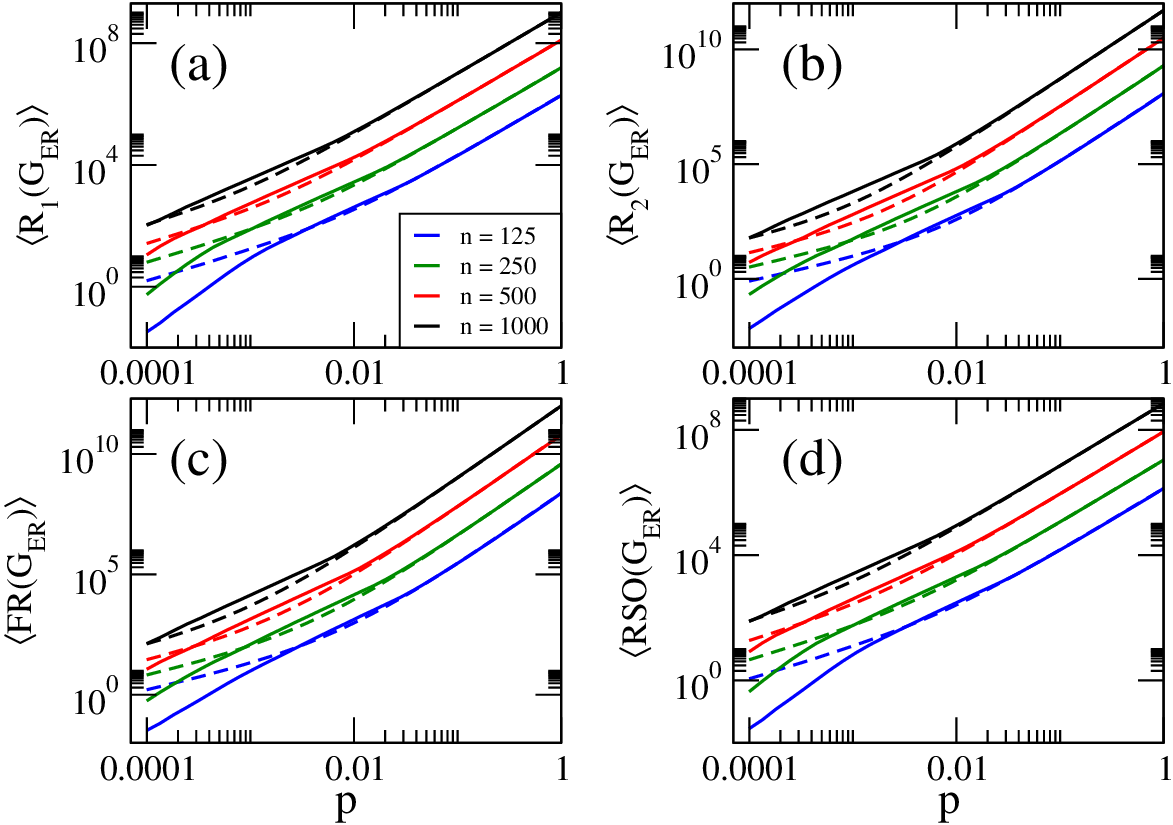}
\caption{(a) $\bra R_1(G_{\tbox{ER}}) \ket$, (b) $\bra R_2(G_{\tbox{ER}}) \ket$, 
(c) $\bra FR(G_{\tbox{ER}}) \ket$, and (d) $\bra RSO(G_{\tbox{ER}}) \ket$
as a function of the probability $p$ of Erd\H{o}s-R\'{e}nyi graphs $G_{\tbox{ER}}(n,p)$ of sizes 
$n\in[125,1000]$.
Dashed lines are (a) $\bra M_1(G_{\tbox{ER}}) \ket$, (b) $\bra M_2(G_{\tbox{ER}}) \ket$, 
(c) $\bra F(G_{\tbox{ER}}) \ket$, and (d) $\bra SO(G_{\tbox{ER}}) \ket$.
Each data value was computed by averaging over $10^{6}$ random graphs $G_{\tbox{ER}}(n,p)$.}
\label{Fig01}
\end{figure}

In Fig.~\ref{Fig01} we present the average values of the Revan-degree--based TIs
$R_1(G_{\tbox{ER}})$, $R_2(G_{\tbox{ER}})$, $FR(G_{\tbox{ER}})$ and $RSO(G_{\tbox{ER}})$
as a function of the probability $p$ of ER graphs of four different sizes $n$ (full lines).
For comparison purposes in each panel of Fig.~\ref{Fig01} we include the corresponding average 
degree--based TIs; that is, we plot the average values of $M_1(G_{\tbox{ER}})$, $M_2(G_{\tbox{ER}})$, 
$F(G_{\tbox{ER}})$ and $SO(G_{\tbox{ER}})$, respectively (dashed lines).

\begin{figure}
\centering
\includegraphics[width=0.8\textwidth]{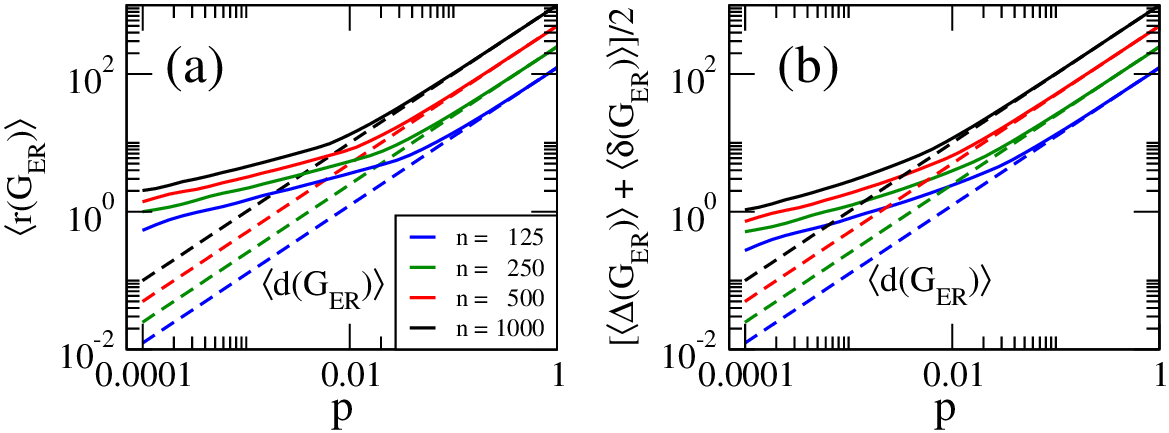}
\caption{(a) Average Revan vertex degree $\bra r(G_{\tbox{ER}}) \ket$ and 
(b) $[\bra \Delta(G_{\tbox{ER}}) \ket + \bra \delta(G_{\tbox{ER}}) \ket]/2$
as a function of the probability $p$ of Erd\H{o}s-R\'{e}nyi graphs $G_{\tbox{ER}}(n,p)$ of 
sizes $n\in[125,1000]$.
Dashed lines in (a,b) are the corresponding average degrees $\bra d(G_{\tbox{ER}}) \ket$.
Each data value was computed by averaging over $10^{6}$ random graphs $G_{\tbox{ER}}(n,p)$.}
\label{Fig02}
\end{figure}

It is interesting to note, from Fig.~\ref{Fig01}, that 
$\bra RX_\Sigma(G_{\tbox{ER}}) \ket \approx \bra X_\Sigma(G_{\tbox{ER}}) \ket$ once $p> 0.01$.
Moreover, given that $RX_\Sigma(G_{\tbox{ER}})$ and $X_\Sigma(G_{\tbox{ER}})$ have the same 
functional form on $r$ and $d$,
respectively, $\bra RX_\Sigma(G_{\tbox{ER}}) \ket \approx \bra X_\Sigma(G_{\tbox{ER}}) \ket$ must be the 
consequence of
\begin{equation}
\bra r(G_{\tbox{ER}}) \ket = \bra \Delta(G_{\tbox{ER}}) \ket + \bra \delta(G_{\tbox{ER}}) \ket - \bra d(G_{\tbox{ER}}) \ket 
\approx \bra d(G_{\tbox{ER}}) \ket ,
\label{rd}
\end{equation}
for large $p$.
Indeed, in Fig.~\ref{Fig02}(a) we plot $\bra r(G_{\tbox{ER}}) \ket$ (full lines) and $\bra d(G_{\tbox{ER}}) \ket$
(dashed lines) and clearly verify that $\bra r(G_{\tbox{ER}}) \ket \approx \bra d(G_{\tbox{ER}}) \ket$ for large $p$.
Thus, the approximation in Eq.~(\ref{rd}) implies that 
$\bra d(G_{\tbox{ER}}) \ket \approx [\bra \Delta(G_{\tbox{ER}}) \ket + \bra \delta(G_{\tbox{ER}}) \ket]/2$.
This rough estimate of the mean from the max and min values is validated in Fig.~\ref{Fig02}(b) where
we contrast $[\bra \Delta(G_{\tbox{ER}}) \ket + \bra \delta(G_{\tbox{ER}}) \ket]/2$ with $\bra d(G_{\tbox{ER}}) \ket$
and show that they certainly coincide for large enough $p$.

Therefore, in the dense limit, i.e.~when $\bra d(G_{\tbox{ER}}) \ket \gg 1$, we can estimate the 
Revan-degree--based TIs by the use of the approximations $r_u \approx r_v \approx \bra r(G_{\tbox{ER}}) \ket$
and $\bra d(G_{\tbox{ER}}) \ket \approx \bra r(G_{\tbox{ER}}) \ket$. For example, for $R_1(G_{\tbox{ER}})$
we can write
$$
R_1(G_{\tbox{ER}}) = \sum_{uv\in E(G_{\tbox{ER}})} r_u + r_v 
\approx \sum_{uv\in E(G_{\tbox{ER}})} 2\bra r(G_{\tbox{ER}}) \ket
\approx n \bra d(G_{\tbox{ER}}) \ket \bra r(G_{\tbox{ER}}) \ket
\approx n \bra r(G_{\tbox{ER}}) \ket^2
$$
or
\begin{equation}
\frac{R_1(G_{\tbox{ER}})}{n} \approx \bra r(G_{\tbox{ER}}) \ket^2 .
\label{R1}
\end{equation}
Above we have used $|E(G_{\tbox{ER}})|=n\bra d(G_{\tbox{ER}}) \ket/2$.
Similar approximations give
\begin{equation}
\frac{R_2(G_{\tbox{ER}})}{n} \approx \frac{1}{2} \bra r(G_{\tbox{ER}}) \ket^3 ,
\label{R2}
\end{equation}
\begin{equation}
\frac{FR(G_{\tbox{ER}})}{n} \approx \bra r(G_{\tbox{ER}}) \ket^3 ,
\label{FR}
\end{equation}
and
\begin{equation}
\frac{RSO(G_{\tbox{ER}})}{n} \approx \frac{1}{\sqrt{2}} \bra r(G_{\tbox{ER}}) \ket^2 .
\label{RSO}
\end{equation}
From Eqs.~(\ref{R1}-\ref{RSO}) we can see that the ratio $RX_\Sigma(G_{\tbox{ER}})/n$ should
depend on $\bra r(G_{\tbox{ER}}) \ket$ only in the dense limit.

\begin{figure}
\centering
\includegraphics[width=0.8\textwidth]{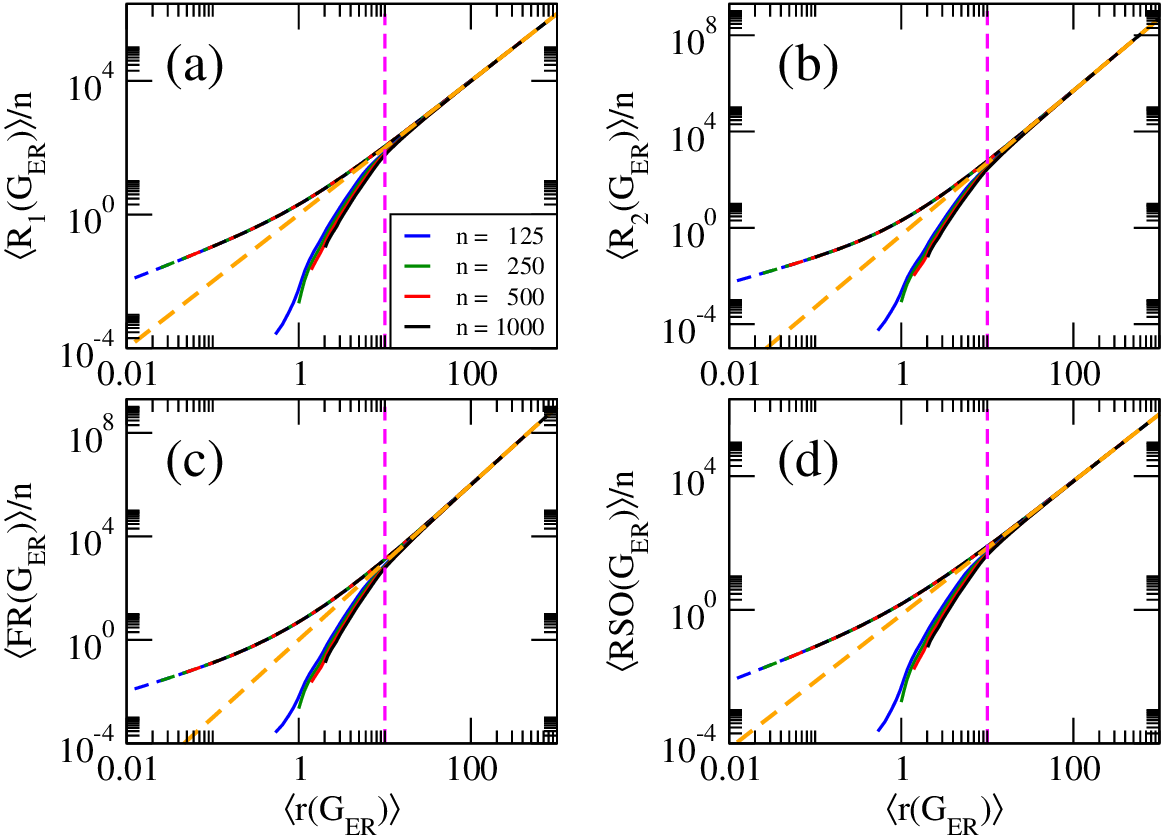}
\caption{(a) $\bra R_1(G_{\tbox{ER}}) \ket/n$, (b) $\bra R_2(G_{\tbox{ER}}) \ket/n$, 
(c) $\bra FR(G_{\tbox{ER}}) \ket/n$, and (d) $\bra RSO(G_{\tbox{ER}}) \ket/n$
as a function of the average Revan vertex degree $\bra r(G_{\tbox{ER}}) \ket$ of Erd\H{o}s-R\'{e}nyi 
graphs $G_{\tbox{ER}}(n,p)$ of sizes $n\in[125,1000]$.
Dashed lines are (a) $\bra M_1(G_{\tbox{ER}}) \ket/n$, (b) $\bra M_2(G_{\tbox{ER}}) \ket/n$, 
(c) $\bra F(G_{\tbox{ER}}) \ket/n$, and (d) $\bra SO(G_{\tbox{ER}}) \ket/n$
as a function of the average degree $\bra d(G_{\tbox{ER}}) \ket$.
Same data of Fig.~\ref{Fig01}.
Orange dashed lines are (a) Eq.~(\ref{R1}), (b) Eq.~(\ref{R2}), (c) Eq.~(\ref{FR}), and (d) Eq.~(\ref{RSO}).
The vertical magenta dashed lines indicate $\bra r(G_{\tbox{ER}}) \ket=10$.}
\label{Fig03}
\end{figure}

Then, in Fig.~\ref{Fig03} we plot $\left< RX_\Sigma(G_{\tbox{ER}}) \right>/n$ vs.~$\left< r(G_{\tbox{ER}}) \right>$ 
(full lines) and observe a good correspondence with Eqs.~(\ref{R1}-\ref{RSO}) (orange dashed lines) in the
dense limit, i.e.~when $\bra r(G_{\tbox{ER}}) \ket \ge 10$.
Furthermore, except for a small-size effect evident at small $\left< r(G_{\tbox{ER}}) \right>$, we notice that the 
curves $\left< RX_\Sigma(G_{\tbox{ER}}) \right>/n$ vs.~$\left< r(G_{\tbox{ER}}) \right>$ do not depend on $n$ 
(that is, the curves for different graph sizes fall one on top of the other) even for $\bra r(G_{\tbox{ER}}) \ket < 10$.
Therefore, a scaling relation for $\left< RX_\Sigma(G_{\tbox{ER}}) \right>$ can be stated as 
\begin{equation}
\frac{\left< RX_\Sigma(G_{\tbox{ER}}) \right>}{n} \approx g_\Sigma(\left< r(G_{\tbox{ER}}) \right>) \, .
\label{scaling2}
\end{equation}

Note that scaling (\ref{scaling2}) is the Revan version of scaling (\ref{scaling1}). Also note that in those
expressions we deliberately named the functions on the rhs as $g_\Sigma$ and $f_\Sigma$, respectively, to 
stress that they are different. Nevertheless, as can be clearly seen in Fig.~\ref{Fig03} where we also include the 
curves $\left< X_\Sigma(G_{\tbox{ER}}) \right>/n$ vs.~$\left< d(G_{\tbox{ER}}) \right>$ (dashed lines), 
once $\bra r(G_{\tbox{ER}}) \ket \ge 10$ the curves $\left< X_\Sigma(G_{\tbox{ER}}) \right>/n$ 
vs.~$\left< d(G_{\tbox{ER}}) \right>$ and $\left< RX_\Sigma(G_{\tbox{ER}}) \right>/n$ vs.~$\left< r(G_{\tbox{ER}}) \right>$ 
coincide.
This means that Eqs.~(\ref{R1}-\ref{RSO}) with $RX\to X$ and $r\to d$ also describe the corresponding 
degree--based indices $X_\Sigma(G_{\tbox{ER}})$ when $\bra d(G_{\tbox{ER}}) \ket \ge 10$; or equivalently,
the functions $f_\Sigma$ and $g_\Sigma$ in the scalings (\ref{scaling1}) and (\ref{scaling2}), respectively, must 
be equal in the dense limit.

\subsection{Revan-degree--based TIs on random geometric graphs}
\label{Sec:TIonER}

\begin{figure}
\centering
\includegraphics[width=0.8\textwidth]{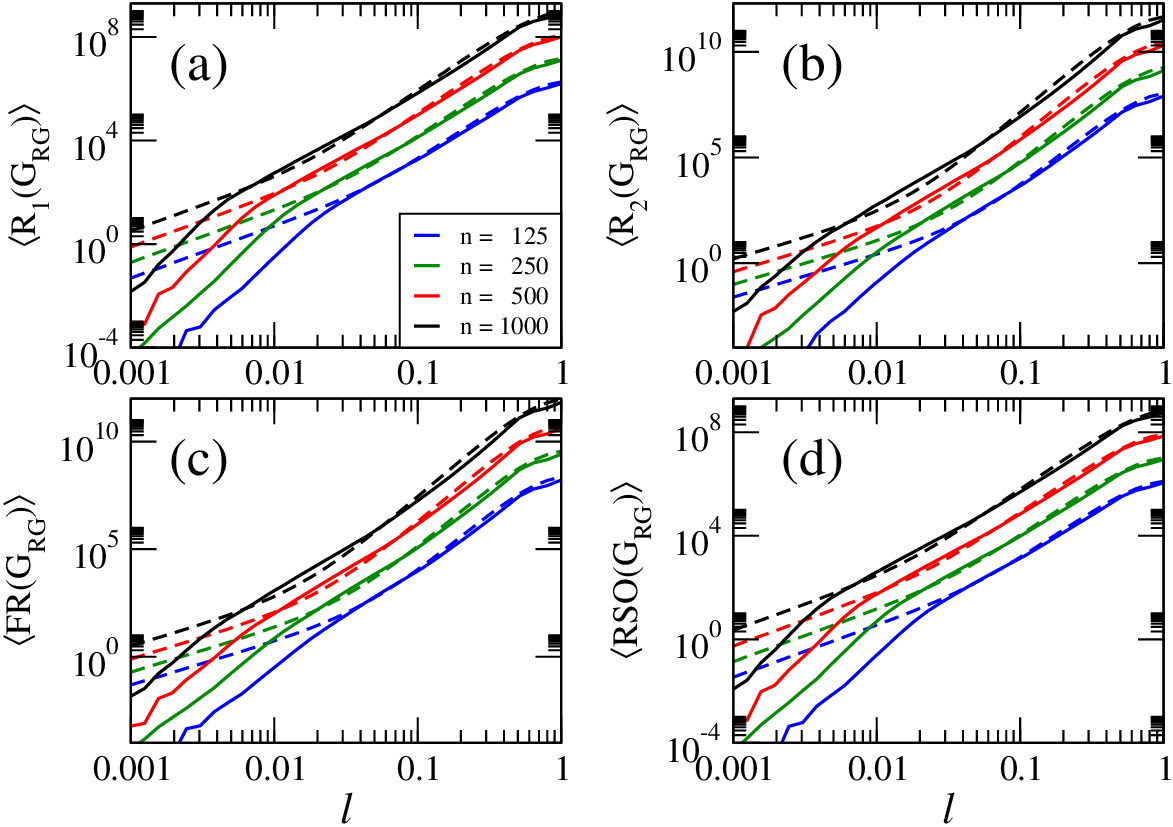}
\caption{(a) $\bra R_1(G_{\tbox{RG}}) \ket$, (b) $\bra R_2(G_{\tbox{RG}}) \ket$, 
(c) $\bra FR(G_{\tbox{RG}}) \ket$, and (d) $\bra RSO(G_{\tbox{RG}}) \ket$
as a function of the connection radius $\ell$ of random geometric graphs $G_{\tbox{RG}}(n,\ell)$ of sizes 
$n\in[125,1000]$.
Dashed lines are (a) $\bra M_1(G_{\tbox{RG}}) \ket$, (b) $\bra M_2(G_{\tbox{RG}}) \ket$, 
(c) $\bra F(G_{\tbox{RG}}) \ket$, and (d) $\bra SO(G_{\tbox{RG}}) \ket$.
Each data value was computed by averaging over $10^{6}$ random graphs $G_{\tbox{RG}}(n,\ell)$.}
\label{Fig04}
\end{figure}

\begin{figure}
\centering
\includegraphics[width=0.8\textwidth]{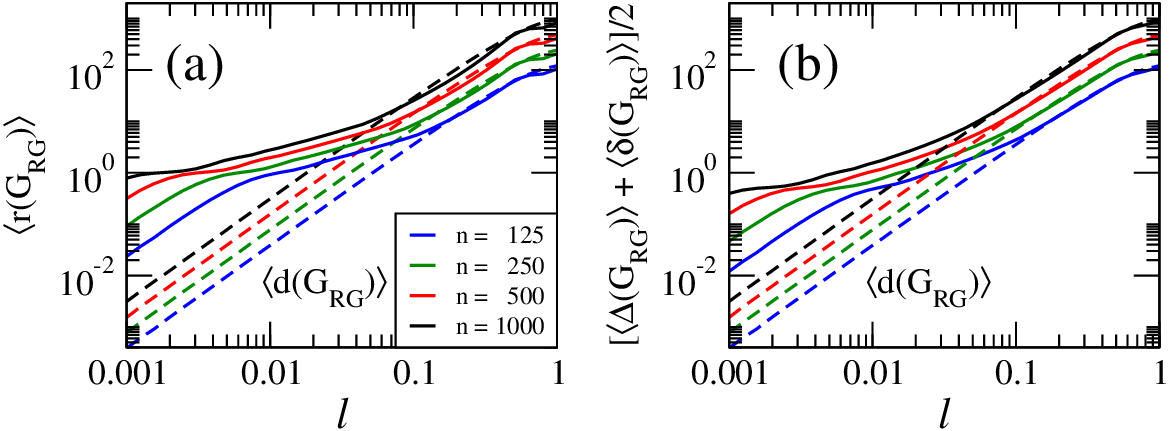}
\caption{(a) Average Revan vertex degree $\bra r(G_{\tbox{RG}}) \ket$ and  
(b) $[\bra \Delta(G_{\tbox{RG}}) \ket + \bra \delta(G_{\tbox{RG}}) \ket]/2$
 as a function of the connection radius $\ell$ of random geometric graphs $G_{\tbox{RG}}(n,\ell)$ of 
sizes $n\in[125,1000]$.
Dashed lines are the corresponding average degrees $\bra d(G_{\tbox{RG}}) \ket$.
Each data value was computed by averaging over $10^{6}$ random graphs $G_{\tbox{RG}}(n,\ell)$.}
\label{Fig05}
\end{figure}

Now, in Fig.~\ref{Fig04} we present the average values of the Revan-degree--based TIs
$R_1(G_{\tbox{RG}})$, $R_2(G_{\tbox{RG}})$, $FR(G_{\tbox{RG}})$ and $RSO(G_{\tbox{RG}})$
as a function of the connection radius $\ell$ of RG graphs of four different sizes $n$ (full lines).
In this figure we also include the corresponding average degree--based TIs as dashed lines.
In addition, in Figs.~\ref{Fig05}(a) and~\ref{Fig05}(b) we plot $\bra r(G_{\tbox{RG}}) \ket$ (full lines) and 
$[\bra \Delta(G_{\tbox{RG}}) \ket + \bra \delta(G_{\tbox{RG}}) \ket]/2$ (full lines) as a function of $\ell$, 
respectively.

\begin{figure}
\centering
\includegraphics[width=0.8\textwidth]{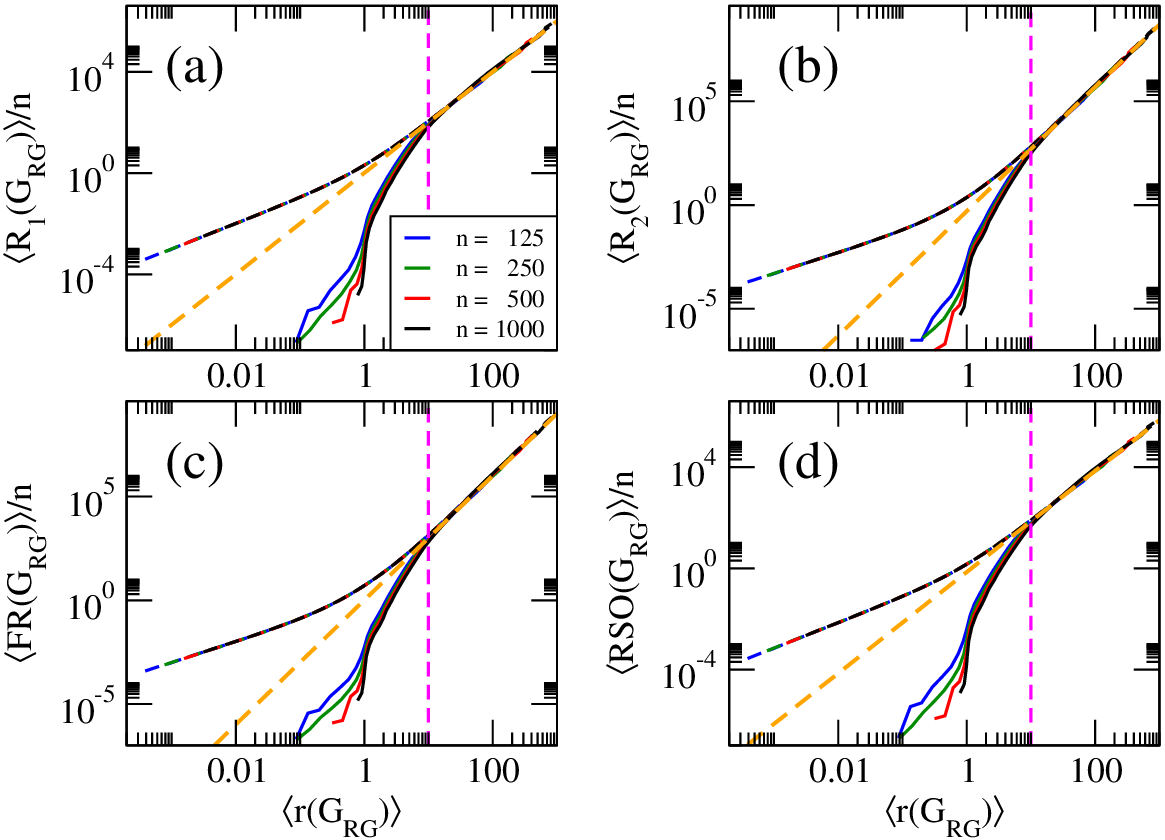}
\caption{(a) $\bra R_1(G_{\tbox{RG}}) \ket/n$, (b) $\bra R_2(G_{\tbox{RG}}) \ket/n$, 
(c) $\bra FR(G_{\tbox{RG}}) \ket/n$, and (d) $\bra RSO(G_{\tbox{RG}}) \ket/n$
as a function of the average Revan vertex degree $\bra r(G_{\tbox{RG}}) \ket$ of random geometric graphs 
$G_{\tbox{RG}}(n,\ell)$ of sizes $n\in[125,1000]$.
Dashed lines are (a) $\bra M_1(G_{\tbox{RG}}) \ket/n$, (b) $\bra M_2(G_{\tbox{RG}}) \ket/n$, 
(c) $\bra F(G_{\tbox{RG}}) \ket/n$, and (d) $\bra SO(G_{\tbox{RG}}) \ket/n$
as a function of the average degree $\bra d(G_{\tbox{RG}}) \ket$.
Same data of Fig.~\ref{Fig04}.
Orange dashed lines are (a) Eq.~(\ref{R1}), (b) Eq.~(\ref{R2}), (c) Eq.~(\ref{FR}), and (d) Eq.~(\ref{RSO})
with $G_{\tbox{ER}} \to G_{\tbox{RG}}$.
The vertical magenta dashed lines indicate $\bra r(G_{\tbox{RG}}) \ket=10$.}
\label{Fig06}
\end{figure}

For comparison purposes, Figs.~\ref{Fig01} and~\ref{Fig02} for ER graphs are equivalent to 
Figs.~\ref{Fig04} and~\ref{Fig05} for RG graphs. In fact, all observations and conclusions made in the previous 
Subsection for ER graphs are also valid for RG graphs, namely:
\begin{itemize}

\item[(i)]
$\bra RX_\Sigma(G_{\tbox{RG}}) \ket \approx \bra X_\Sigma(G_{\tbox{RG}}) \ket$ for large $\ell$, see Fig.~\ref{Fig04},

\item[(ii)]
$\bra d(G_{\tbox{RG}}) \ket \approx [\bra \Delta(G_{\tbox{RG}}) \ket + \bra \delta(G_{\tbox{RG}}) \ket]/2$
for large $\ell$, see Fig.~\ref{Fig05}(b), thus

\item[(iii)]
$\bra r(G_{\tbox{RG}}) \ket \approx \bra d(G_{\tbox{RG}}) \ket$ for large $\ell$, see Fig.~\ref{Fig05}(a). Therefore,

\item[(iv)]
Eqs.~(\ref{R1}-\ref{RSO}) with $G_{\tbox{ER}}\to G_{\tbox{RG}}$ should also be valid for RG graphs in the dense
limit. This is indeed verified in Fig.~\ref{Fig06} where we have plotted 
$\left< RX_\Sigma(G_{\tbox{RG}}) \right>/n$ vs.~$\left< r(G_{\tbox{RG}}) \right>$ (full lines) 
together with Eqs.~(\ref{R1}-\ref{RSO}) (orange dashed lines).

\item[(v)] Finally, from Fig.~\ref{Fig06}, the scaling law  
\begin{equation}
\frac{\left< RX_\Sigma(G_{\tbox{RG}}) \right>}{n} \approx g_\Sigma(\left< r(G_{\tbox{RG}}) \right>)
\label{scaling4}
\end{equation}
can be stated. 

\end{itemize}

\section{Statistical analysis of Revan-degree--based MTIs on random graphs}
\label{Sec:MTI}
 
We now introduce the multiplicative versions of the Revan-degree--based TIs, $RX_\Pi(G)$, studied in the
previous Section:
$$
R_{1\Pi}(G) = \prod_{uv\in E(G)} r_u + r_v \, ,
\qquad
R_{2\Pi}(G) = \prod_{uv\in E(G)} r_u r_v \, ,
$$
$$
FR_\Pi(G) = \prod_{uv\in E(G)} r_u^2 + r_v^2 \, ,
$$
and
$$
RSO_\Pi(G) = \prod_{uv\in E(G)} \sqrt{r_u^2 + r_v^2} \, .
$$
These MTIs are the Revan versions of the multiplicative Zagreb indices~\cite{EIG12,K16}
$$
\Pi_1^*(G) = \prod_{uv\in E(G)} d_u + d_v \, ,
\qquad
\Pi_2(G) = \prod_{uv\in E(G)} d_u d_v \, ,
$$
the multiplicative forgotten index
\begin{equation}
F_\Pi(G) = \prod_{uv\in E(G)} d_u^2 + d_v^2 \, ,
\label{multF}
\end{equation}
and the multiplicative Sombor index
\begin{equation}
SO_\Pi(G) = \prod_{uv\in E(G)} \sqrt{d_u^2 + d_v^2} \, ,
\label{multSO}
\end{equation}
respectively.
Note that (as far as we know) neither $F_\Pi(G)$ nor $SO_\Pi(G)$ have been considered before.

It is fair to recall that a statistical study of degree--based MTIs, $X_\Pi(G)$, on random graphs 
has been already reported in Ref.~\cite{AMMS22}. There, the multiplicative Zagreb indices, the 
multiplicative Randi\'c connectivity index, the multiplicative harmonic index, the multiplicative 
sum-connectivity index, the multiplicative inverse degree index, as well as the Narumi-Katayama index 
were applied to ER graphs, RG graphs and bipartite random graphs. There it was demonstrated that
$\left< \ln X_\Pi(G) \right>$ normalized to the order of the graph scale with the corresponding average 
degree:
\begin{equation}
\frac{\left< \ln X_\Pi(G) \right>}{n} \equiv f_\Pi(\left< d \right>) \ .
\label{scaling5}
\end{equation}
Note that scaling (\ref{scaling5}) can be considered as the multiplicative version of scaling (\ref{scaling1}).

Then, in what follows we apply the Revan-degree--based MTIs defined above on ER and RG graphs.

\begin{figure}
\centering
\includegraphics[width=0.8\textwidth]{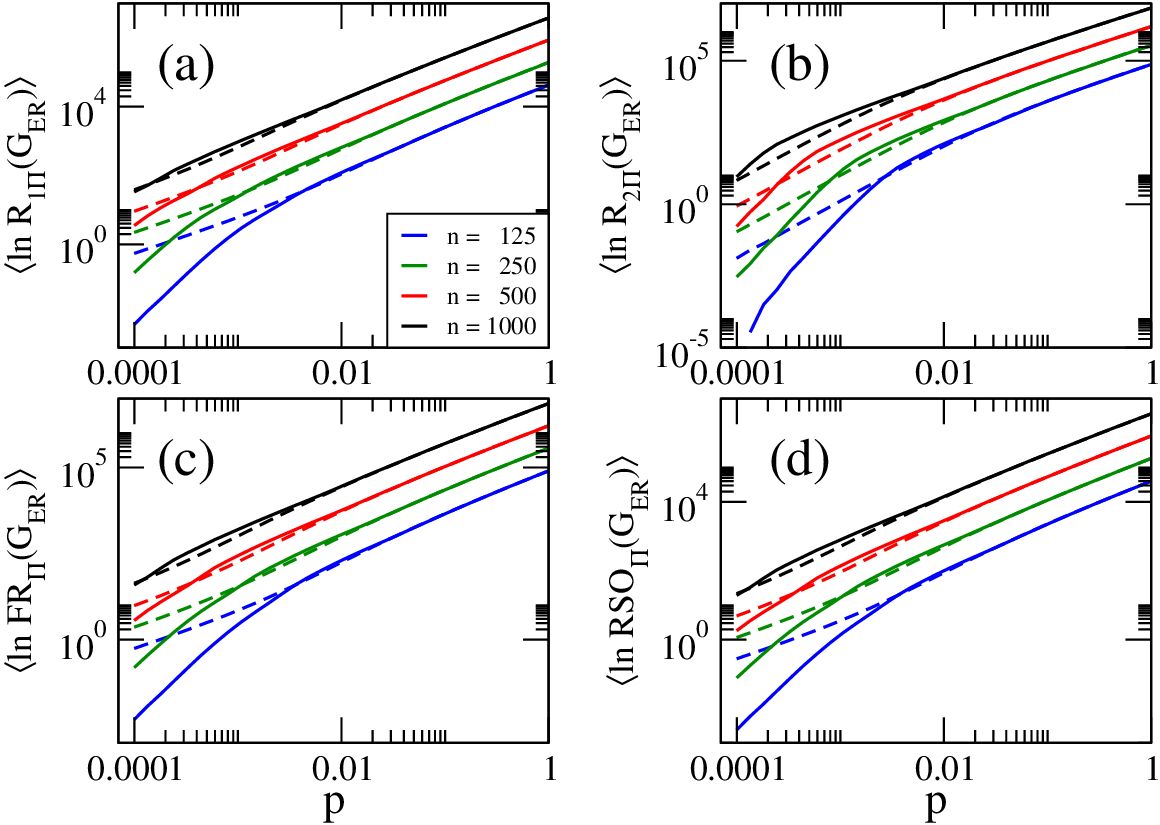}
\caption{(a) $\bra \ln R_{1\Pi}(G_{\tbox{ER}}) \ket$, (b) $\bra \ln R_{2\Pi(}G_{\tbox{ER}}) \ket$, 
(c) $\bra \ln FR_\Pi(G_{\tbox{ER}}) \ket$, and (d) $\bra \ln RSO_\Pi(G_{\tbox{ER}}) \ket$
as a function of the probability $p$ of Erd\H{o}s-R\'{e}nyi graphs $G_{\tbox{ER}}(n,p)$ of sizes 
$n\in[125,1000]$.
Dashed lines are (a) $\bra \ln \Pi^*_1(G_{\tbox{ER}}) \ket$, (b) $\bra \ln \Pi_2(G_{\tbox{ER}}) \ket$, 
(c) $\bra \ln F_\Pi(G_{\tbox{ER}}) \ket$, and (d) $\bra \ln SO_\Pi(G_{\tbox{ER}}) \ket$.
Each data value was computed by averaging over $10^{6}$ random graphs $G_{\tbox{ER}}(n,p)$.}
\label{Fig07}
\end{figure}

\begin{figure}
\centering
\includegraphics[width=0.8\textwidth]{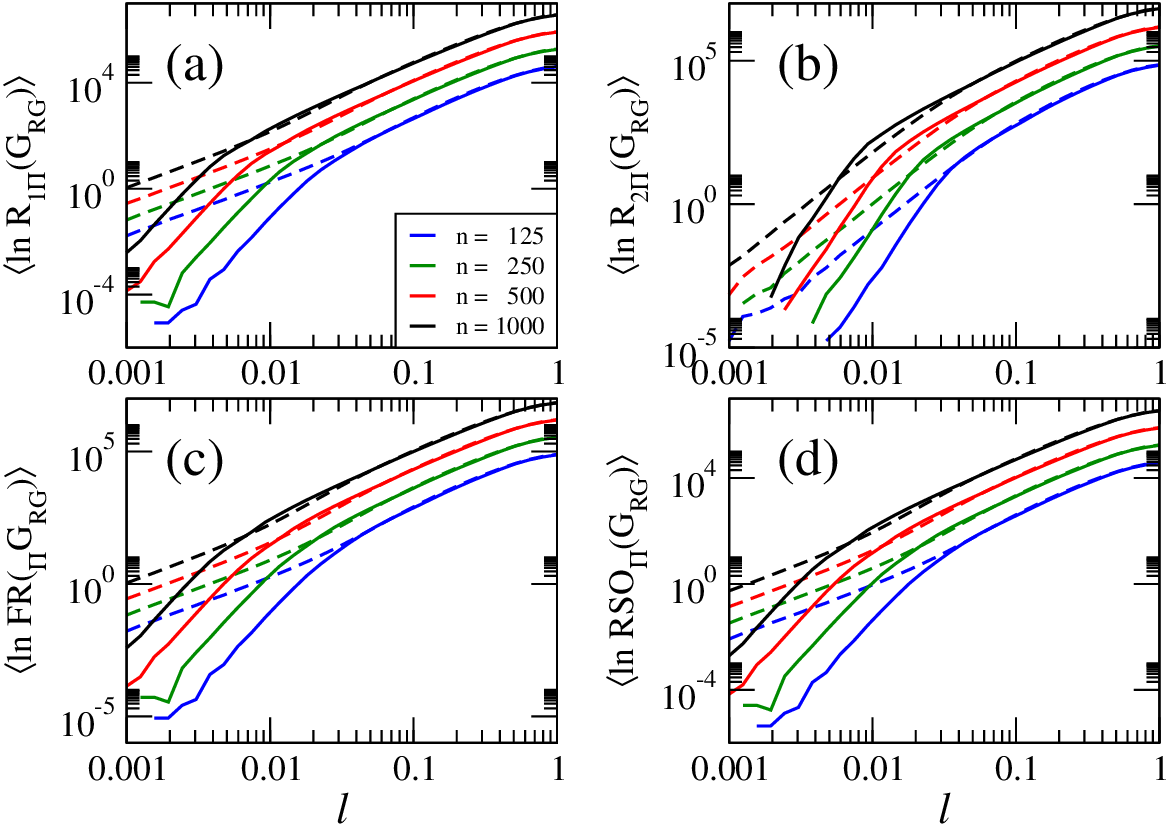}
\caption{(a) $\bra \ln R_{1\Pi}(G_{\tbox{RG}}) \ket$, (b) $\bra \ln R_{2\Pi}(G_{\tbox{RG}}) \ket$, 
(c) $\bra \ln FR_\Pi(G_{\tbox{RG}}) \ket$, and (d) $\bra \ln RSO_\Pi(G_{\tbox{RG}}) \ket$
as a function of the connection radius $\ell$ of random geometric graphs $G_{\tbox{RG}}(n,\ell)$ of sizes 
$n\in[125,1000]$.
Dashed lines are (a) $\bra \ln \Pi^*_1(G_{\tbox{RG}}) \ket$, (b) $\bra \ln \Pi_2(G_{\tbox{RG}}) \ket$, 
(c) $\bra \ln F_\Pi(G_{\tbox{RG}}) \ket$, and (d) $\bra \ln SO_\Pi(G_{\tbox{RG}}) \ket$.
Each data value was computed by averaging over $10^{6}$ random graphs $G_{\tbox{RG}}(n,\ell)$.}
\label{Fig09}
\end{figure}

In Figs.~\ref{Fig07} and~\ref{Fig09} we present the average values of the Revan-degree--based 
MTIs $R_{1\Pi}(G)$, $R_{2\Pi}(G)$, $FR_\Pi(G)$ and $RSO_\Pi(G)$ for ER (Fig.~\ref{Fig07}) and RG 
(Fig.~\ref{Fig09}) graphs of four different sizes $n$ (full lines).
In these figures we also plot the corresponding average degree--based MTIs as dashed lines.
For both random graph models we observe that $\bra RX_\Pi(G) \ket \approx \bra X_\Pi(G) \ket$ 
for large enough $p$ and large enough $\ell$, respectively.


Indeed, as for the Revan-degree--based TIs, here we can also estimate the Revan-degree--based 
MTIs in the dense limit by the use of the approximations $r_u \approx r_v \approx \bra r(G) \ket$
and $\bra d(G) \ket \approx \bra r(G) \ket$. Thus, for $R_{1\Pi}(G)$ we write
$$
R_{1\Pi}(G) = \prod_{uv\in E(G)} r_u + r_v 
\approx \prod_{uv\in E(G)} 2\bra r(G) \ket
\approx [2\bra r(G) \ket]^{n\bra d(G) \ket/2}
$$
which leads to
$$
\ln R_{1\Pi}(G) 
\approx \frac{1}{2} n\bra d(G) \ket \ln [2\bra r(G) \ket]
\approx \frac{1}{2} n\bra r(G) \ket \ln [2\bra r(G) \ket]
$$
or
\begin{equation}
\frac{\ln R_{1\Pi}(G)}{n} \approx \frac{1}{2} \bra r(G) \ket \ln [2\bra r(G) \ket] .
\label{RR1}
\end{equation}
Similar approximations give
\begin{equation}
\frac{\ln R_{2\Pi}(G)}{n} \approx \bra r(G) \ket \ln \bra r(G) \ket ,
\label{RR2}
\end{equation}
\begin{equation}
\frac{\ln FR_\Pi(G)}{n} \approx \bra r(G) \ket \ln [\sqrt{2}\bra r(G) \ket] ,
\label{RFR}
\end{equation}
and
\begin{equation}
\frac{\ln RSO_\Pi(G)}{n} \approx \frac{1}{2} \bra r(G) \ket \ln [\sqrt{2}\bra r(G) \ket] .
\label{RRSO}
\end{equation}
Note that Eqs.~(\ref{RR1}-\ref{RRSO}) should work for both ER and RG graphs.

\begin{figure}
\centering
\includegraphics[width=0.8\textwidth]{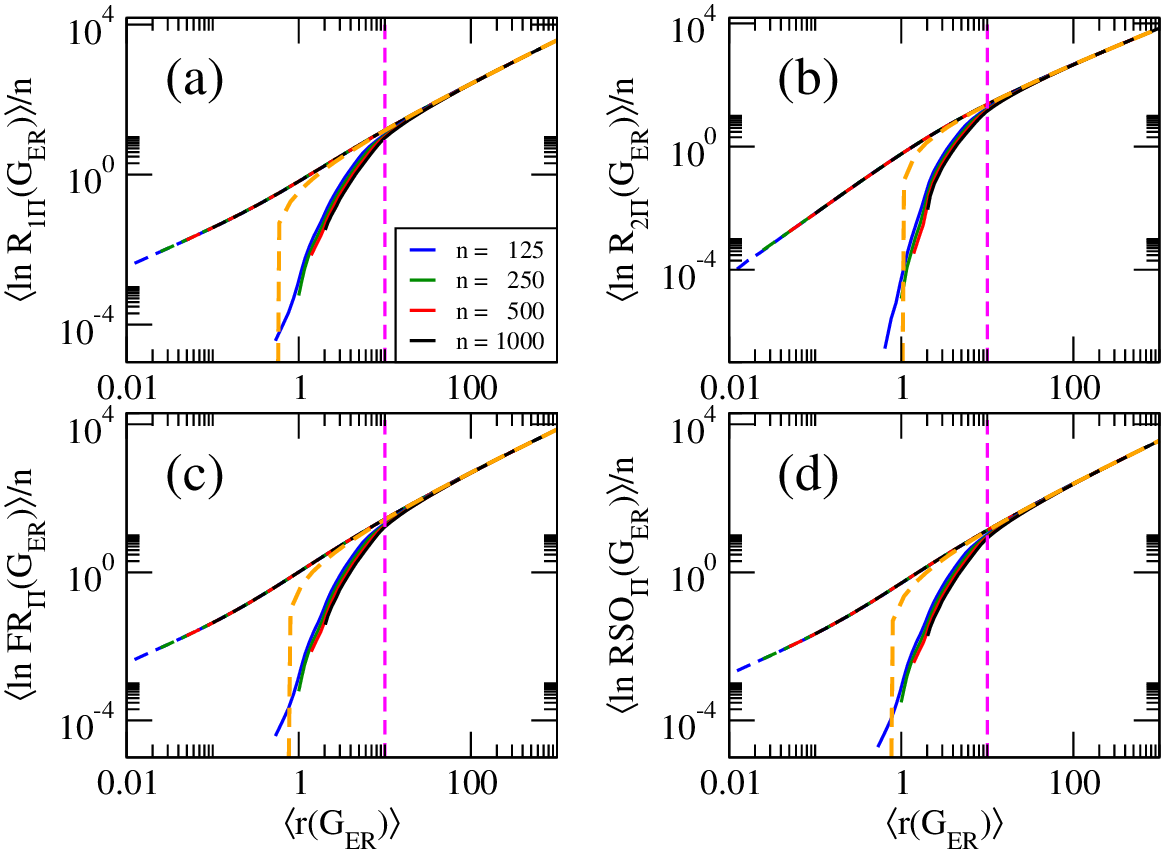}
\caption{(a) $\bra \ln R_{1\Pi}(G_{\tbox{ER}}) \ket/n$, (b) $\bra \ln R_{2\Pi}(G_{\tbox{ER}}) \ket/n$, 
(c) $\bra \ln FR_\Pi(G_{\tbox{ER}}) \ket/n$, and (d) $\bra \ln RSO_\Pi(G_{\tbox{ER}}) \ket/n$
as a function of the average Revan vertex degree $\bra r(G_{\tbox{ER}}) \ket$ of Erd\H{o}s-R\'{e}nyi 
graphs $G_{\tbox{ER}}(n,p)$ of sizes $n\in[125,1000]$.
Dashed lines are (a) $\bra \ln \Pi^*_1(G_{\tbox{ER}}) \ket/n$, (b) $\bra \ln \Pi_2(G_{\tbox{ER}}) \ket/n$, 
(c) $\bra \ln F_\Pi(G_{\tbox{ER}}) \ket/n$, and (d) $\bra \ln SO_\Pi(G_{\tbox{ER}}) \ket/n$
as a function of the average degree $\bra d(G_{\tbox{ER}}) \ket$.
Same data of Fig.~\ref{Fig07}.
Orange dashed lines are (a) Eq.~(\ref{RR1}), (b) Eq.~(\ref{RR2}), (c) Eq.~(\ref{RFR}), and (d) Eq.~(\ref{RRSO}).
The vertical magenta dashed lines indicate $\bra r(G_{\tbox{ER}}) \ket=10$.}
\label{Fig08}
\end{figure}
\begin{figure}
\centering
\includegraphics[width=0.8\textwidth]{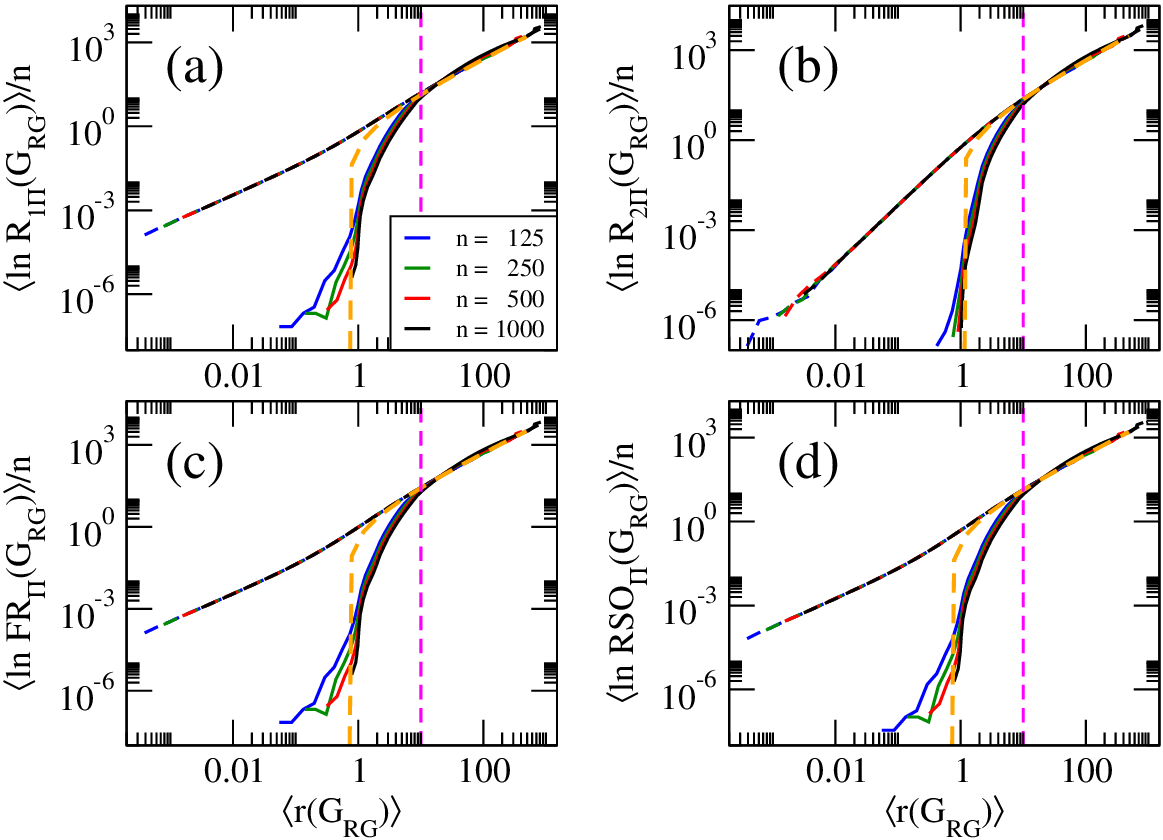}
\caption{(a) $\bra \ln R_{1\Pi}(G_{\tbox{RG}}) \ket/n$, (b) $\bra \ln R_{2\Pi}(G_{\tbox{RG}}) \ket/n$, 
(c) $\bra \ln FR_\Pi(G_{\tbox{RG}}) \ket/n$, and (d) $\bra \ln RSO_\Pi(G_{\tbox{RG}}) \ket/n$
as a function of the average Revan vertex degree $\bra r(G_{\tbox{RG}}) \ket$ of random geometric graphs 
$G_{\tbox{RG}}(n,\ell)$ of sizes $n\in[125,1000]$.
Dashed lines are (a) $\bra \ln \Pi^*_1(G_{\tbox{RG}}) \ket/n$, (b) $\bra \ln \Pi_2(G_{\tbox{RG}}) \ket/n$, 
(c) $\bra \ln F_\Pi(G_{\tbox{RG}}) \ket/n$, and (d) $\bra \ln SO_\Pi(G_{\tbox{RG}}) \ket/n$
as a function of the average degree $\bra d(G_{\tbox{RG}}) \ket$.
Same data of Fig.~\ref{Fig09}.
Orange dashed lines are (a) Eq.~(\ref{RR1}), (b) Eq.~(\ref{RR2}), (c) Eq.~(\ref{RFR}), and (d) Eq.~(\ref{RRSO}).
The vertical magenta dashed lines indicate $\bra r(G_{\tbox{RG}}) \ket=10$.}
\label{Fig10}
\end{figure}

Therefore, in Figs.~\ref{Fig08} and~\ref{Fig10} we plot $\left< RX_\Pi(G) \right>/n$ vs.~$\left< r(G) \right>$ 
(full lines) for ER and RG graphs, respectively, together with Eqs.~(\ref{RR1}-\ref{RRSO}) (orange dashed lines). 
Indeed, we observe a very good correspondence between predictions~(\ref{RR1}-\ref{RRSO}) and the numerical 
data for both random graph models in dense limit, i.e.~when $\bra r(G) \ket \ge 10$.
From these figures, except for a small-size effect at small $\left< r(G_{\tbox{ER}}) \right>$, we can state
the scaling of the Revan-degree--based MTIs as
\begin{equation}
\frac{\left< \ln RX_\Pi(G) \right>}{n} \approx g_\Pi(\left< r(G) \right>) \, .
\label{scaling6}
\end{equation}
Note that scalings~(\ref{scaling5}) and~(\ref{scaling6}) indeed coincide for $\bra r(G) \ket \ge 10$
as can be clearly seen in Figs.~\ref{Fig08} and~\ref{Fig10} where we also include the 
curves $\left< X_\Pi(G_{\tbox{ER}}) \right>/n$ vs.~$\left< d(G_{\tbox{ER}}) \right>$ (dashed lines).
This means that Eqs.~(\ref{RR1}-\ref{RRSO}) with $RX\to X$ and $r\to d$ also describe the corresponding 
degree--based indices $X_\Pi(G)$ when $\bra d(G) \ket \ge 10$; or equivalently,
the functions $f_\Pi$ and $g_\Pi$ in the scalings (\ref{scaling5}) and (\ref{scaling6}), respectively, must be equal in 
the dense limit.

\section{Summary}

Motivated by potential theoretical--practical applications of topological indices, in this work we perform 
a thorough numerical study of two families of Revan-degree--based graph invariants: 
$RX_\Sigma(G) = \sum_{uv \in E(G)} F(r_u,r_v)$ and 
$RX_\Pi(G) = \prod_{uv \in E(G)} F(r_u,r_v)$. In particular while $RX_\Sigma(G)$ has gained interest 
recently, see e.g.~Refs.~\cite{K17,K18,KG22,KMRS22}, we are introducing $RX_\Pi(G)$ here.
Specifically, we have considered the Revan-degree--based versions of the first and second Zagrev 
indices, the forgotten index, and the Sombor index.
We have applied both $X_\Sigma(G)$ and $X_\Pi(G)$ on ensembles of Erd\"os-R\'enyi graphs and 
random geometric graphs, see Figs.~\ref{Fig01}, \ref{Fig04}, \ref{Fig07} and~\ref{Fig09}.

We would like to add that we have also introduced here the multiplicative forgotten index 
$F_\Pi(G)$ and the multiplicative Sombor index $SO_\Pi(G)$, see Eqs.~(\ref{multF}) and~(\ref{multSO}), 
respectively.

On the one hand we have shown that $\left< RX_\Sigma(G) \right>$ and $\left< \ln RX_\Pi(G) \right>$, 
normalized to the order of the graph $n$, scale with the average Revan degree $\left< r \right>$; that is,
\begin{equation}
\frac{\left< RX_\Sigma(G) \right>}{n} \approx g_\Sigma(\left< r(G) \right>) \qquad \mbox{and} \qquad
\frac{\left< \ln RX_\Pi(G) \right>}{n} \approx g_\Pi(\left< r(G) \right>) \, ,
\label{scalings}
\end{equation}
see Figs.~\ref{Fig03}, \ref{Fig06}, \ref{Fig08} and~\ref{Fig10}.
On the the other hand we have provided expressions for both $\left< RX_\Sigma(G) \right>$ and 
$\left< \ln RX_\Pi(G) \right>$
in the dense graph limit, see Eqs.~(\ref{R1}-\ref{RSO}) and Eqs.~(\ref{RR1}-\ref{RRSO}), respectively.

In addition, we have found that $\bra r(G) \ket \approx \bra d(G) \ket$ and
$\bra RX_\Sigma(G) \ket \approx \bra X_\Sigma(G) \ket$ in the dense limit, i.~e.~when $\bra r(G) \ket \ge 10$.
This makes the scalings in (\ref{scalings}) to reproduce the scalings (\ref{scaling1}) and (\ref{scaling5}) reported
in Refs.~\cite{MMRS20,AHMS20,MMRS21,AMRS21,MAAS22} and~\cite{AMMS22}, respectively. 
Therefore, Eqs.~(\ref{R1}-\ref{RSO}) and Eqs.~(\ref{RR1}-\ref{RRSO}) also describe the corresponding 
degree-based topological indices in the dense limit.
Furthermore, it is relevant to stress that the clear difference between Revan-degree--based indices and
degree--based indices for $\bra r(G) \ket < 10$, makes Revan-degree--based indices particularly useful 
in that regime where they could provide additional information to standard degree--based indices.

We hope that our study may motivate further computational as well as theoretical studies of 
Revan-degree--based topological indices.

\section*{\bf ACKNOWLEDGEMENTS}
J.A.M.-B. thanks support from CONACyT (Grant No. 286633) and VIEP-BUAP (Grant No. 100405811-VIEP2022), Mexico.
J.M.S. acknowledges financial support from
Agencia Estatal de Investigaci\'on (PID2019-106433GB-I00 / AEI / 10.13039/501100011033), Spain.

%\section*{\bf REFERENCES}

\end{document}